# On the Well-Posedness of Green's Function Reconstruction via the Kirchhoff-Helmholtz Equation for One-Speed Neutron Diffusion

Roberto Ponciroli[1]

Nuclear Science and Engineering Division, Argonne National Laboratory,
9700 S. Cass Avenue, Lemont, IL, 60439, USA

## ABSTRACT

This work presents a methodology for reconstructing the spatial distribution of the neutron flux in a nuclear reactor, leveraging real-time measurements obtained from ex-core detectors. The Kirchhoff-Helmholtz (K-H) equation inherently defines the problem of estimating a scalar field within a domain based on boundary data, making it a natural mathematical framework for this task. The main challenge lies in deriving the Green's function specific to the domain and the neutron diffusion process. While analytical solutions for Green's functions exist for simplified geometries, their derivation of complex, heterogeneous domains—such as a nuclear reactor—requires a numerical approach.

The objective of this work is to demonstrate the well-posedness of the data-driven Green's function approximation by formulating and solving the K-H equation as an inverse problem. After establishing the symmetry properties that the Green's function must satisfy, the K-H equation is derived from the one-speed neutron diffusion model. This is followed by a comprehensive description of the procedure for interpreting sensor readings and implementing the neutron flux reconstruction algorithm. Finally, the existence and uniqueness of the Green's function inferred from the sampled data are demonstrated, ensuring the reliability of the proposed method and its predictions.

## 1. INTRODUCTION

Traditionally, nuclear reactors employ hundreds of in-core detectors to monitor neutron flux [1]. However, in advanced reactor designs, placing sensors within the core is often impractical due to spatial constraints and extreme operating conditions. As a result, these systems require alternative approaches capable of reconstructing the neutron flux distribution in real-time using a limited number of sensors, preferably located at the periphery of the reactor. The objective of this work is to develop an algorithm capable of estimating the spatial distribution of the neutron flux ($\phi(\vec{x})$) in real-time using data from ex-core detectors. To this aim, it is necessary to estimate the Green's function that describes the propagation of neutrons within a nuclear reactor.

The fundamental assumption is that the one-speed neutron diffusion model provides a valid framework for interpreting both simulation results and experimental measurements of the neutron flux distribution governed by the neutron transport equation [2]. Specifically, neutron energy dependence is neglected by assuming that all neutrons can be characterized by a single kinetic energy. This assumption is justified since ex-core detectors in use cannot discriminate neutron energies, so the added precision of a multi-group model cannot be leveraged given the available measurement limitations. In addition, the linearity of the neutron diffusion equation enables the application of the superposition principle, exploits inherent symmetries, and allows the use of well-established mathematical techniques for its solution. Even so, computing the Green's function continues to pose significant challenges. A nuclear reactor is a heterogeneous domain, consisting of regions such as fuel, cladding, coolant, and graphite, each exhibiting significantly different macroscopic cross-section properties. Although neutron diffusion principles remain valid in all regions, the Green's function cannot be derived analytically due to the geometric complexity and the medium's spatial heterogeneity. Numerical methods represent the only viable approach, and machine learning techniques can be particularly effective when guided by a first-principles formulation. The proposed methodology is designed for application to open bounded domains in $\mathbb{R}^d$, where $d = 2,3$. In this study, we denote the monitored domain by $\Omega$, and its enclosing boundary—typically a surface—by $\partial\Omega \equiv \Gamma$. The objective is to establish a mapping that takes any two points $\vec{x}, \vec{y} \in \Omega$ and returns the corresponding Green's function value $\hat{G}(\vec{x}, \vec{y})$, as specified in Eq.(1).

---

[1] Correspondence: rponciroli@anl.gov; Tel.: +1-630-252-3455

$$(\vec{x}, \vec{y}) \xrightarrow{\text{mapping}} \hat{G}(\vec{x}, \vec{y}) \tag{1}$$

To accomplish this, the Kirchhoff-Helmholtz (K-H) integral equation for the one-speed neutron diffusion model is employed. Traditionally, the K-H equation is used to solve boundary value problems (BVPs), wherein the scalar field at any point within the domain is determined based on the Green's function associated with the diffusion operator and the prescribed Dirichlet and Neumann boundary conditions (BCs) [3]. In this work, the K-H equation is utilized to address the inverse problem—reconstructing the corresponding Green's function based on the measured BCs and the field values at a finite set of points within the domain.

In Section 2, the one-speed neutron diffusion equation is presented, and the corresponding source term characteristic of control rod maneuvering is defined. In Section 3, it is demonstrated that Dirichlet problems with various non-homogeneous BCs can always be reformulated as an equivalent homogeneous problem using a procedure known as 'lifting'. This transformation is crucial in demonstrating the versatility of the proposed approach. Additionally, the self-adjointness of the associated operator is demonstrated. This property ensures specific symmetry characteristics of the corresponding Green's function, which are exploited in the numerical algorithm for its data-driven estimation. In Section 4, the K-H equation from the one-speed neutron diffusion model is derived, along with the steps of the procedure for interpreting sensor readings and reconstructing the neutron flux distribution. In Section 5, the existence and uniqueness of the Green's function inferred from the sampled data in the inverse problem is demonstrated. In Section 6, the validity of the results obtained for a constant source term was extended to the nonlinear case, in which the source is proportional to the local neutron flux—thereby demonstrating the method's applicability under realistic reactor conditions. Finally, the main conclusion are drawn (Section 7).

## 2. DESCRIPTION OF THE ONE-SPEED NEUTRON DIFFUSION MODEL

In Eq.(2), the one-speed neutron diffusion equation for steady-state conditions (no external neutron source is present) is reported. Eqs.(3)(4) define the macroscopic transport cross section ($\Sigma_{tr}(\vec{y})$) and the diffusion coefficient ($D(\vec{y})$), respectively. The parameters used in diffusion equation are listed in Table 1.

$$\nabla \cdot (D(\vec{y})\nabla\phi(\vec{y})) - \Sigma_a(\vec{y})\phi(\vec{y}) - \Sigma_a^{CR}(\vec{y})\phi(\vec{y}) + \nu\Sigma_f(\vec{y})\phi(\vec{y}) = 0 \quad \vec{y} \in \Omega \tag{2}$$

$$\Sigma_{tr}(\vec{y}) = \Sigma_t(\vec{y}) - \bar{\mu}_0\Sigma_s(\vec{y}) \tag{3}$$

$$D(\vec{y}) = \frac{1}{3\Sigma_{tr}(\vec{y})} \tag{4}$$

Across different scenarios, the position of the control rods is the only adjustable parameter influencing the system's response. The fission term cannot be directly controlled, as the concentration of fissile material decreases due to fuel burnup - a process that is beyond operator intervention in conventional reactors. The spatial distribution of the fission cross-section ($\Sigma_f(\vec{y})$) is a fixed property of the medium through which neutrons diffuse (Eq.(5)). The expression of the source term in our applications ($Q(\vec{y})$) is provided in Eq.(6). It comprises a fixed component (the fission contribution, $\nu\Sigma_f(\vec{y})$) and a controllable component (the absorption cross-section in the region containing the control rods, $\delta\Sigma_a^{CR}(\vec{y})$, as expressed in Eq.(7)).

$$-\nabla \cdot \left(D(\vec{y})\nabla\phi(\vec{y})\right) + \Sigma_a(\vec{y})\phi(\vec{y}) = \nu\Sigma_f(\vec{y})\phi(\vec{y}) - \delta\Sigma_a^{CR}(\vec{y})\phi(\vec{y}) \qquad (5)$$

$$Q = \left(\nu\Sigma_f(\vec{y}) - \delta\Sigma_a^{CR}(\vec{y})\right)\phi(\vec{y}) \qquad (6)$$

$$\delta\Sigma_a^{CR}(\vec{y}) = \begin{cases} \Sigma_a^{CR} & \text{for } \vec{y} \in \text{Rod} \\ 0 & \text{for } \vec{y} \notin \text{Rod} \end{cases} \qquad (7)$$

The intensity of the adopted source term is proportional to the reconstructed neutron flux in the corresponding regions ($\phi(\vec{y})$). Given the mathematical complexity introduced by the nonlinear source term, we begin by analyzing a simplified case in which the source is modeled as a spatial function independent of the neutron flux ($Q = Q(\vec{y}) \in L^2(\Omega)$). Once this case is validated, the results will be extended to a more realistic scenario (Section 6).

**Table 1. List of parameters used in the neutron diffusion model (Eqs.(2)(3)(4)).**

| Parameter | Definition |
|---|---|
| $D(\vec{y})$ | Neutron diffusion coefficient (m) |
| $\Sigma_a(\vec{y})$ | Macroscopic absorption cross section of structural and fuel materials (m$^{-1}$) |
| $\Sigma_a^{CR}(\vec{y})$ | Macroscopic absorption cross-section due to control rods (m$^{-1}$) |
| $\nu$ | Average number of neutrons emitted per fission event (-) |
| $\Sigma_f(\vec{y})$ | Macroscopic fission cross-section (m$^{-1}$) |
| $\Sigma_{tr}(\vec{y})$ | Macroscopic transport cross-section (m$^{-1}$) |
| $\Sigma_t(\vec{y})$ | Total macroscopic cross-section (m$^{-1}$) |
| $\Sigma_s(\vec{y})$ | Macroscopic scattering cross-section (m$^{-1}$) |
| $\bar{\mu}_0$ | Average cosine of the scattering angle (-) |

By substituting Eq.(6) into Eq.(5), the non-homogeneous diffusion equation is obtained (Eq.(8)). For simplicity, the spatial dependence of the variables is omitted.

$$-\nabla \cdot (D\nabla\phi) + \Sigma_a\phi = Q \qquad (8)$$

The neutron diffusion operator for the current application ($\mathcal{L}$) is reported in Eq.(9).

$$\mathcal{L}\phi = -\nabla \cdot (D\nabla\phi) + \Sigma_a\phi \qquad (9)$$

The Green's function $G(\vec{x}, \vec{y})$ for the linear differential operator $\mathcal{L}$ is defined as the fundamental solution satisfying Eq.(10). It represents the response at $\vec{x}$ to a unit impulse applied at $\vec{y}$.

$$-\nabla \cdot (D\nabla G) + \Sigma_a G = \delta(\vec{x} - \vec{y}) \qquad (10)$$

The Green's function is generally defined as a function of two variables, i.e., the field point and the source point.

- **Field Point**: it is the first argument ($\vec{x}$). This is the point in the domain where the effect of a source (or BCs) is being observed or measured.
- **Boundary Variable** (or **Source Point**): it is the second argument ($\vec{y}$). This variable represents the location of the source or the integration variable over the boundary. It is often integrated over when constructing the solution to ensure that the solution satisfies the prescribed BCs.

In the context of reconstructing the neutron flux from sensor readings, the source points refer to locations in the peripheral region of the core—or even outside of it—where ex-core detectors are positioned. The field points are locations inside the core where direct measurement is not feasible, but where the neutron flux needs to be estimated. A fundamental property of the Green's function is that it depends solely on the BCs and the intrinsic dynamics of the system, remaining independent of any external forcing term. This aspect will play a crucial role in ensuring the versatility of the flux reconstruction method proposed in this work. Further details on this topic will be provided in Section 3.2.

## 3. TREATMENT OF BOUNDARY CONDITIONS IN THE STUDIED MONITORING PROBLEM

### 3.1. Role of boundary conditions on neutron flux reconstruction within the reactor core

The proposed methodology aims to develop a Digital Twin (DT) for neutron flux reconstruction within a nuclear reactor. This term refers to an algorithm that integrates a parameterized representation of the system's physics, enabling predictive capabilities for both control and diagnostics [4]. During the training phase, the DT parameters are optimized to characterize the system's reference conditions; in the testing phase, the trained algorithm is supplied with real-time measurements to generate predictions. The algorithm combines limited in-core sensor readings with an extensive array of ex-core measurements, so it is critical that data collected at the reactor's exterior boundary carry enough information to infer the interior flux distribution. This aspect has two major implications:

- The neutron flux and its gradient will never be zero (*non-homogeneous BCs*).
- BCs vary across different operational scenarios. Variations in control rod positions and power levels lead to different neutron flux distributions, which in turn produce distinct boundary values at the domain's outer surface. This variability is advantageous, as it allows identifying different operational conditions using ex-core measurements. If the BCs were constant, measuring the neutron flux or its gradient at the boundaries would be meaningless, as they would provide no relevant information about the internal reactor conditions.

To quantitatively illustrate these concepts, we examine the results of a scoping study conducted in MCNP6 [5] to compute the neutron flux distribution within the Purdue University Reactor Number One (PUR-1) core (Figure 1). A comprehensive description of the reactor model and its technical specifications can be found in [6]. To highlight the importance of measurements conducted in the peripheral region, the neutron flux profile along the x-axis of the reactor core was sampled. Virtual sensors were placed in the gap between two adjacent fuel assemblies, as illustrated by the thin vertical wires in Figure 2a. The simulated neutron flux detected by each virtual sensor is shown in Figure 2b. Results indicate that sensors positioned within the reflector region measure approximately 40% of the neutron flux found at the core center, while those located on the outer surface of the reflector record about 30% of the core-center flux. These findings demonstrate that the reflector does not function as a complete neutron shield for the PUR-1 core. Despite the expected attenuation, a substantial portion of the thermal neutron flux remains detectable outside the core boundary. This enables the reconstruction of the neutron flux distribution within the active core region, making these ex-core measurements viable for core monitoring applications.

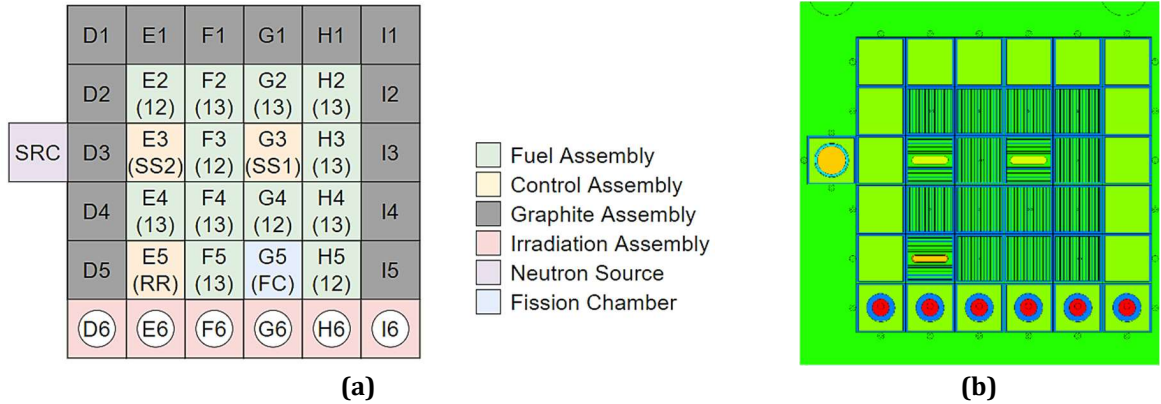

**Figure 1.** (a) PUR-1 core layout. (b) Graphical representation of the geometry adopted in the MCNP6 model of PUR-1 (ex-core detectors are displayed) [6].

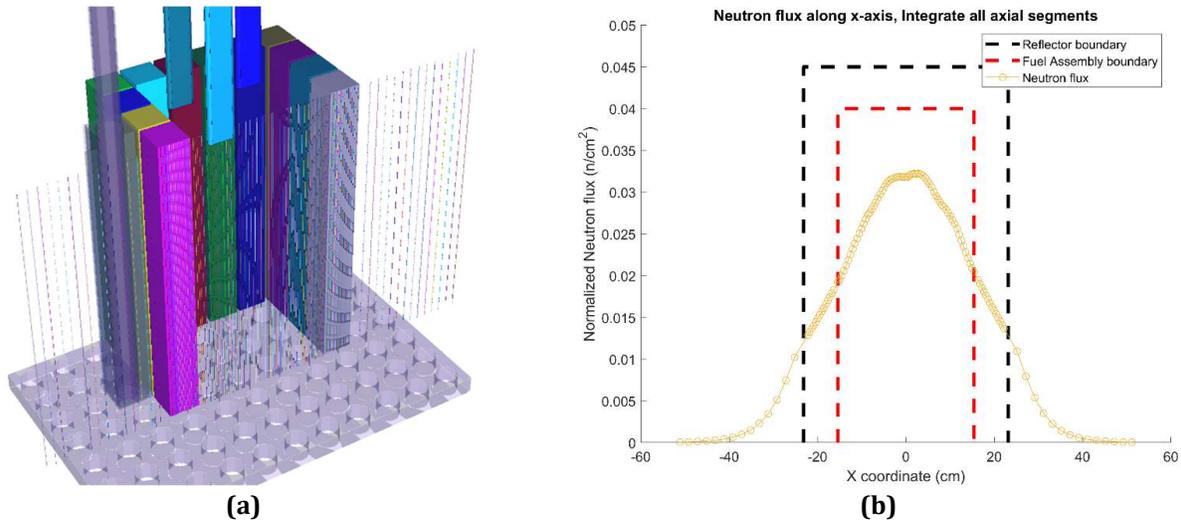

**Figure 2.** (a) Placement of virtual sensors along the x-axis in the PUR-1 core model. (b) Simulated neutron flux distribution along the x-axis of PUR-1 [6].

*3.2. Approach for managing boundary condition variability across different operational scenarios*

BCs play a fundamental role in ensuring the uniqueness of the Green's function. They define how the Green's function satisfies the governing physical constraints of the system, thereby ensuring a physically consistent and well-posed solution [7]. As stressed in Section 3.1, in our applications, the BCs are non-homogeneous and vary across different operational scenarios. Specifically, each change in control rod positions induces a corresponding modification in the neutron flux distribution throughout the reactor, along with a variation in both the neutron flux and its gradient at the boundary. While this variability is advantageous from a monitoring perspective, it introduces a mathematical challenge. In principle, because the Green's function depends on both the governing operator and the imposed BCs, a unique Green's function would be required for each scenario. To ensure consistent application of the reconstructed Green's function, we propose transforming each non-homogeneous BVP into to a reference homogeneous BVP, which serves as a common framework for all scenarios. The Green's function estimated for this case can then be used to reconstruct the neutron flux across all other scenarios. This transformation technique is known as the *lifting method* [3]. It consists of choosing a function that extends the prescribed BCs into the interior of the domain. A typical choice for $w$ (when the geometry permits) is a linear or polynomial function that matches the boundary data inferred from sensor measurements. The lifting method can be applied to both Dirichlet and Neumann BCs.

- **Measurement Perspective**: from the standpoint of data acquisition, the two transformations are equivalent - both require measuring the neutron flux and its gradient at the boundary.
- **Mathematical Perspective:** applying the lifting method to Neumann conditions introduces complexity. The lifting function must satisfy stricter regularity and compatibility conditions, and the resulting solution is only defined up to an additive constant, reflecting the fact that Neumann problems are inherently non-unique without further normalization. Given these challenges, in this work, the lifting method is applied to Dirichlet BCs.

### 3.3. Proof of well-posedness of non-homogenous Dirichlet problems and applicability of the lifting method

In this Section, we demonstrate the well-posedness of both the non-homogeneous Dirichlet problem (Eq.(11)) and its transformed counterpart obtained via the lifting method.

$$\begin{cases} \mathcal{L}\phi = Q & \text{in } \Omega \\ \phi = \psi & \text{on } \partial\Omega \end{cases} \tag{11}$$

By performing a lifting, the non-homogeneous Dirichlet problem can be transformed into an equivalent one with homogeneous conditions. First, let us define a suitable function ($w$) (Eq.(12)).

$$\begin{aligned} & w \in H^1(\Omega) \\ & w = \psi \text{ on } \partial\Omega \end{aligned} \tag{12}$$

This formulation enables solving for $u$ in the space $H_0^1(\Omega)$, where functions vanish on the boundary $\partial\Omega$. Then, reconstructing the full solution $\phi$ using Eq.(13).

$$\phi = u + w \tag{13}$$

Since the operator $\mathcal{L}$ is linear, substituting Eq.(13) into the original equation and rearranging the terms lead to the formulation of the corresponding homogenous problem (Eq.(14)).

$$\begin{cases} \mathcal{L}u = Q - \mathcal{L}w & \text{in } \Omega \\ u = 0 & \text{on } \partial\Omega \end{cases} \tag{14}$$

Let us derive the variational formulation for $u$. After substituting Eq.(13) into the diffusion equation, we multiply by a test function $v \in H_0^1(\Omega)$ and integrate over $\Omega$ (Eq.(15)).

$$\int_\Omega [-\nabla \cdot (D\nabla(u + w)) + \Sigma_a(u + w)]v \, d\Omega = \int_\Omega Qv \, d\Omega \tag{15}$$

After integrating by parts, we obtain Eq.(16).

$$\int_\Omega D\nabla(u + w) \cdot \nabla v \, d\Omega + \int_\Omega \Sigma_a(u + w)v \, d\Omega = \int_\Omega Qv \, d\Omega \tag{16}$$

This expression naturally splits into terms involving $u$ and those involving $w$. A bilinear form on $H_0^1(\Omega)$ (Eq.(17)) and a linear functional (Eq.(18)) are identified. We seek a solution $u$ for the problem formulated in Eq.(19).

$$a(u, v) = \int_\Omega D\nabla u \cdot \nabla v \, d\Omega + \int_\Omega \Sigma_a uv \, d\Omega \tag{17}$$

$$F(v) = \int_\Omega Qv \, d\Omega - \int_\Omega D\nabla w \cdot \nabla v \, d\Omega - \int_\Omega \Sigma_a wv \, d\Omega \tag{18}$$

$$a(u, v) = F(v) \quad \forall v \in H_0^1(\Omega) \tag{19}$$

The Lax-Milgram lemma guarantees the existence and uniqueness of the solution, provided that the bilinear form $a(u, v)$ is continuous and coercive, and the linear functional $F(v)$ is continuous on $H_0^1(\Omega)$. As for the physical properties of the diffusive medium, we assume that $D(\vec{y})$ is bounded and strictly positive and $\Sigma_a(\vec{y})$ is bounded, i.e., there exist $D_{min} > 0$, $D_{max} > 0$ and $\Sigma_{a,min} \geq 0$, $\Sigma_{a,max} > 0$ such that the conditions in Eq.(20) hold.

$$0 < D_{min} \leq D(\vec{y}) \leq D_{max}, \quad 0 \leq \Sigma_{a,min} \leq \Sigma_a(\vec{y}) \leq \Sigma_{a,max} \quad \vec{y} \in \Omega \tag{20}$$

- *Continuity (Boundedness) of $a(u, v)$*

    We want to demonstrate that there exists a constant $M_a > 0$ such that the condition in Eq.(21) is satisfied.

    $$|a(u, v)| \leq M_a \|u\|_{H_0^1(\Omega)} \|v\|_{H_0^1(\Omega)} \quad \forall u, v \in H_0^1(\Omega) \tag{21}$$

    For the first term, by applying the Cauchy-Schwarz inequality, we obtain Eq.(22).

    $$\left| \int_\Omega D\nabla u \cdot \nabla v \, d\Omega \right| \leq \int_\Omega |D\nabla u \cdot \nabla v| \, d\Omega \leq D_{max} \int_\Omega |\nabla u||\nabla v| \, d\Omega \leq D_{max} \|\nabla u\|_{L^2(\Omega)} \|\nabla v\|_{L^2(\Omega)} \tag{22}$$

    The norm in space $H^1$ accounts for both the $L^2$-norm of the function and its gradient (Eq.(23)).

    $$\|u\|_{H^1(\Omega)} = \left( \|u\|_{L^2(\Omega)}^2 + \|\nabla u\|_{L^2(\Omega)}^2 \right)^{1/2} \tag{23}$$

    For functions in $H_0^1(\Omega)$, the Poincaré inequality holds, ensuring the existence of a constant $C_P > 0$, which depends on $\Omega$, such that the condition in Eq.(24) is satisfied.

    $$\|u\|_{L^2(\Omega)} \leq C_P \|\nabla u\|_{L^2(\Omega)} \tag{24}$$

This inequality shows that the $L^2$-norm of $u$ is controlled by the $L^2$-norm of its gradient. This allows defining a simpler norm on $H_0^1(\Omega)$ (Eq.(25)). Under this norm, the first term is bounded Eq.(26)).

$$\|u\|_{H_0^1(\Omega)} = \|\nabla u\|_{L^2(\Omega)} \tag{25}$$

$$\left|\int_\Omega D\nabla u \cdot \nabla v \, d\Omega\right| \leq D_{max}\|u\|_{H_0^1(\Omega)}\|v\|_{H_0^1(\Omega)} \tag{26}$$

For the second term, by applying the Cauchy-Schwarz inequality, we obtain Eq.(27).

$$\left|\int_\Omega \Sigma_a uv \, d\Omega\right| \leq \Sigma_{a,max}\int_\Omega |uv| \, d\Omega \leq \Sigma_{a,max}\|u\|_{L^2(\Omega)}\|v\|_{L^2(\Omega)} \tag{27}$$

By applying the Poincaré inequality, we demonstrate that the second term is bounded as well (Eq.(28)).

$$\left|\int_\Omega \Sigma_a uv \, d\Omega\right| \leq \Sigma_{a,max}C_P^2\|u\|_{H_0^1(\Omega)}\|v\|_{H_0^1(\Omega)} \tag{28}$$

By combining the two estimates from Eqs.(26)(28), we obtain Eqs.(29)(30) and the continuity of $a(u,v)$ is demonstrated.

$$|a(u,v)| \leq (D_{max} + \Sigma_{a,max}C_P^2)\|u\|_{H_0^1(\Omega)}\|v\|_{H_0^1(\Omega)} = M_a\|u\|_{H_0^1(\Omega)}\|v\|_{H_0^1(\Omega)} \tag{29}$$

$$M_a = D_{max} + \Sigma_{a,max}C_P^2 \tag{30}$$

- *Coercivity of $a(u,v)$*

  We want to demonstrate that there exists a constant $C_a > 0$ such that the condition in Eq.(31) is satisfied.

  $$|a(u,u)| \geq C_a\|u\|_{H_0^1(\Omega)} \quad \forall u \in H_0^1(\Omega) \tag{31}$$

  Using the assumptions on $D(\vec{y})$ and $\Sigma_a(\vec{y})$ expressed in Eq.(20), we can write Eq.(32).

  $$a(u,u) = \int_\Omega D|\nabla u|^2 \, d\Omega + \int_\Omega \Sigma_a|u|^2 \, d\Omega \geq D_{min}\|\nabla u\|_{L^2(\Omega)} + \Sigma_{a,min}\|u\|_{L^2(\Omega)} \tag{32}$$

  Hence, by applying the Poincaré inequality, the coercivity of $a(u,v)$ is demonstrated (Eqs.(33)(34)).

  $$a(u,u) \geq D_{min}\|\nabla u\|_{L^2(\Omega)} + \Sigma_{a,min}\|u\|_{L^2(\Omega)} \geq C_a\|u\|_{H_0^1(\Omega)} \tag{33}$$

  $$C_a = \min\{D_{min}, \Sigma_{a,min}/C_P^2\} \tag{34}$$

- *Continuity of $F(v)$*

  We want to demonstrate that there exists a constant $C_F \geq 0$ such that the condition in Eq.(35) is satisfied.

  $$|F(v)| \leq C_F\|v\|_{H_0^1(\Omega)} \quad \forall v \in H_0^1(\Omega) \tag{35}$$

Applying the Cauchy-Schwarz inequality to Eq.(18) yields Eq.(36).

$$|F(v)| \leq \|Q\|_{L^2(\Omega)}\|v\|_{L^2(\Omega)} + D_{max}\|\nabla w\|_{L^2(\Omega)}\|\nabla v\|_{L^2(\Omega)} + \Sigma_{a,max}\|w\|_{L^2(\Omega)}\|v\|_{L^2(\Omega)} \qquad (36)$$

After applying the Poincaré inequality and collecting the contributions, we obtain Eq.(37) and the continuity of $F(v)$ is demonstrated.

$$|F(v)| \leq \left(\|Q\|_{L^2(\Omega)}C_P + D_{max}\|\nabla w\|_{L^2(\Omega)} + \Sigma_{max}\|w\|_{L^2(\Omega)}C_P\right)\|v\|_{H_0^1(\Omega)} = C_F\|v\|_{H_0^1(\Omega)} \qquad (37)$$

$$C_F = \|Q\|_{L^2(\Omega)}C_P + D_{max}\|\nabla w\|_{L^2(\Omega)} + \Sigma_{max}\|w\|_{L^2(\Omega)}C_P \qquad (38)$$

By constructing the weak formulation in the space $H_0^1(\Omega)$ and verifying the hypotheses of the Lax-Milgram lemma, we established the well-posedness of the non-homogeneous Dirichlet problem for the one-speed neutron diffusion equation. This result also confirms the mathematical soundness of the lifting method for handling the variability of BCs across different scenarios.

### 3.4. Proof of self-adjointness of the operator in the one-speed neutron diffusion equation

In this section and the following one, we will characterize the regularity of the diffusion operator and the symmetry properties of its corresponding Green's function. These properties are derived solely from the physics of the problem, rather than the specific geometry of the analyzed system. The goal is to identify a set of properties that remain valid across all studied scenarios and can be leveraged in the data-driven optimization of the Green's function. This approach accelerates convergence and ensures that the resulting approximation satisfies the problem's intrinsic constraints. In this Section, the self-adjointness of the operator in the one-speed neutron diffusion equation is demonstrated. The condition for $\mathcal{L}$ operator to be self-adjoint (or "symmetric", since we are dealing with real-valued functions) is presented in Eq.(39). It states that $\mathcal{L}$ must be symmetric with respect to the standard $L^2$-inner product (Eq.(40)).

$$\int_\Omega (\mathcal{L}\phi_1)\phi_2 \, d\Omega = \int_\Omega \phi_1(\mathcal{L}\phi_2) \, d\Omega \qquad (39)$$

$$\langle \mathcal{L}\phi_1, \phi_2 \rangle = \langle \phi_1, \mathcal{L}\phi_2 \rangle \qquad (40)$$

for every pair of test functions $\phi_1, \phi_2$ in $\Omega$ satisfying the BC.

When non-homogeneous BCs are imposed, integration by parts yields nonzero boundary contributions. The lifting method is applied to reformulate the problem as a homogeneous Dirichlet problem (Eq.(14)). This transformation absorbs the boundary contributions into the source term, preserving the symmetry of the operator acting on the modified unknown function ($u$). Besides, the diffusion coefficient is discontinuous due to material heterogeneities in typical reactor physics applications. This leads to a decomposition of the domain into multiple subdomains, each characterized by uniform material properties. <u>The goal is to establish that the diffusion operator remains piecewise self-adjoint across regions with uniform properties. This requires verifying that the integral formulation preserves self-adjointness when applied over the entire domain, accounting for all subdomains.</u>

To demonstrate the self-adjointness of $\mathcal{L}$ operator, we apply integration by parts and show that the resulting boundary terms either cancel out or vanish. First, we substitute the definition of $\mathcal{L}$ from Eq.(9) into Eq.(39), omit the interior integrals that cancel out and arrive at Eq.(41).

$$-\int_\Omega \phi_2 \nabla \cdot (D\nabla \phi_1) d\Omega = -\int_\Omega \phi_1 \nabla \cdot (D\nabla \phi_2) d\Omega \tag{41}$$

We now focus on the interior integral on the left-hand side. Applying integration by parts yields Eq.(42).

$$-\int_\Omega \phi_2 \nabla \cdot (D\nabla \phi_1) \, d\Omega = \int_\Omega D\nabla \phi_1 \nabla \phi_2 \, d\Omega - \int_\Gamma D\phi_2 (\nabla \phi_1 \cdot \vec{n}) dS_y - \sum \int_{\Gamma_i} \ldots \tag{42}$$

where $\Gamma_i$ denotes the internal interfaces between subdomains with different values of the diffusion coefficient ($D_i$).

When this procedure is applied within each material subdomain and then sum over all regions, the resulting surface integrals include contributions from both the domain's outer boundary and the internal interfaces separating adjacent subdomains. On the right-hand side of Eq.(42), the first integral cancels out when integrating by parts the right-hand side of Eq.(41), while the second integral vanishes because of homogeneous Dirichlet BC. Let us consider the third term and denote $\Omega_i$ one of the $N$ subdomains. For each subdomain, a boundary integral term of the form given in Eq.(43) is obtained.

$$\int_{\partial\Omega_i} D_i \phi_2 (\nabla \phi_1 \cdot \vec{n}_i) d\partial\Omega_i \tag{43}$$

where $\vec{n}_i$ is the outward normal to $\Omega_i$.

The internal interface $\Gamma_i$ between adjacent subdomains $\Omega_i$ and $\Omega_j$ is a portion of both $\partial\Omega_i$ and $\partial\Omega_j$, and it is distinct from the outer boundary of the domain. When summing over all subdomains, interface terms of the form given in Eq.(44) emerge.

$$-\int_{\Gamma_i} \phi_2 (D_i \nabla \phi_1 \cdot \vec{n}_i) \, d\Gamma_i \quad \text{and} \quad -\int_{\Gamma_i} \phi_2 (D_j \nabla \phi_1 \cdot \vec{n}_j) \, d\Gamma_i \tag{44}$$

where $D_i$ and $D_j$ are the diffusion coefficients on the two sides of the interface, and $\vec{n}_i = -\vec{n}_j$ is the normal pointing outward from each subdomain.

For the interface integrals on both sides of Eq.(41) to cancel each other, continuity conditions must be imposed on the neutron flux ($\phi$) and the neutron current ($D\nabla \phi \cdot \vec{n}$) across $\Gamma_i$. In this way, no net interface term remains.

- **Flux Continuity**: Neutron flux is assumed to be continuous across material interfaces, as an abrupt change would violate the physical conservation of neutrons.
- **Current Continuity**: The neutron current (which depends on the diffusion coefficient and the flux gradient) is also assumed to be continuous, ensuring that there is no artificial accumulation or loss of neutrons at the interface.

These conditions are not only theoretically justified but are also standard assumptions in reactor physics [2]. The elimination of these contributions, ensuring that no net interface term remains, confirms the self-adjointness of the operator in a heterogenous medium.

### 3.5. Symmetry properties of the Green's function of neutron diffusion operator

As mentioned in Section 2, the properties of a Green's function depend on both the underlying differential operator and the BC. If $\mathcal{L}$ is self-adjoint with respect to the inner product defined on the domain (Eq.(40)) and the BCs are chosen so that the associated bilinear form is symmetric, then the Green's function will satisfy the *reciprocity* relation. This means that swapping the source and field points does not change the Green's function, i.e., the response at point $\vec{x}$ due to a unit impulse at point $\vec{y}$ is identical to the response at $\vec{y}$ when the impulse is applied at $\vec{x}$ (Eq.(45)).

$$G(\vec{x}, \vec{y}) = G(\vec{y}, \vec{x}) \tag{45}$$

This is a consequence of the piecewise self-adjointness of the diffusion operator across regions with uniform properties. The Green's function for the operator $\mathcal{L}$ is the solution of Eq.(46) and satisfies homogeneous Dirichlet BC (Eq.(47)).

$$\mathcal{L}G(\vec{x}, \vec{y}) = \delta(\vec{x} - \vec{y}) \quad \vec{x} \in \Omega \tag{46}$$
$$G(\vec{x}, \vec{y}) = 0 \quad \vec{x} \in \partial\Omega \tag{47}$$

To demonstrate reciprocity, let us test the definition with two points $\vec{x}_1$ and $\vec{x}_2$ within the domain $\Omega$. By applying the self-adjointness condition, we arrive at Eq.(48).

$$\langle G(\vec{y}, \vec{x}_1), \mathcal{L}G(\vec{y}, \vec{x}_2) \rangle = \langle \mathcal{L}G(\vec{y}, \vec{x}_1), G(\vec{y}, \vec{x}_2) \rangle \quad \vec{x}_1, \vec{x}_2 \in \Omega \tag{48}$$

By substituting $\mathcal{L}G(\vec{y}, \vec{x}_1) = \delta(\vec{y} - \vec{x}_1)$ and $\mathcal{L}G(\vec{y}, \vec{x}_2) = \delta(\vec{y} - \vec{x}_2)$ into Eq.(48) and performing the integration, we obtain the reciprocity relation expressed in Eq.(49).

$$G(\vec{x}_2, \vec{x}_1) = G(\vec{x}_1, \vec{x}_2) \quad \vec{x}_1, \vec{x}_2 \in \Omega \tag{49}$$

Homogeneous BCs ensure that no extra terms interfere with this symmetry because the contributions from the boundary are zero. Specifically, applying a lifting procedure to transform the one-speed neutron diffusion equation with non-homogeneous Dirichlet BC into an equivalent homogeneous problem ensures the self-adjointness of the operator within the associated homogeneous function space; once the operator is proven to be self-adjoint, the symmetry of the Green's function in the transformed problem follows naturally. The importance of this symmetry in enabling the accurate reconstruction of the Green's function from collected measurements will be examined in detail in Section 5.

## 4. APPLICATION OF K-H EQUATION TO THE FIELD RECONSTRUCTION PROBLEM

### 4.1. Derivation of K-H equation for the one-speed neutron diffusion model

The K-H equation inherently defines the problem of estimating a scalar field within a domain based on boundary data, making it a natural mathematical framework for this reconstruction task. In this Section, the K-H integral equation corresponding to the one-speed neutron diffusion equation is derived. Our analysis focuses on the non-homogeneous Dirichlet problem defined in Eq.(11). The first step involves applying the lifting method to obtain its homogeneous counterpart, as shown in Eq.(14). By expanding the operator $\mathcal{L}$, we obtain Eq.(50).

$$\begin{cases} -\nabla \cdot (D\nabla u) + \Sigma_a u = Q + \nabla \cdot (D\nabla w) - \Sigma_a w & \text{in } \Omega \\ u = 0 & \text{on } \partial\Omega \end{cases} \quad (50)$$

From this equation, we can define the *effective source* ($\tilde{Q}$) (Eq.(51)). Substituting this expression into Eq.(50) reformulates the problem into the form presented in Eq.(52).

$$\tilde{Q} = Q + \nabla \cdot (D\nabla w) - \Sigma_a w \quad \text{in } \Omega \quad (51)$$

$$\begin{cases} -\nabla \cdot (D\nabla u) + \Sigma_a u = \tilde{Q} & \text{in } \Omega \\ u = 0 & \text{on } \partial\Omega \end{cases} \quad (52)$$

Based on the assumptions made about the source term and the parameters of the neutron diffusion equation, we can characterize the properties of the function $\tilde{Q}$.

- $Q \in L^2(\Omega)$ is a spatially varying function that generates and absorbs neutrons independently of the neutron flux.
- $w \in H^1(\Omega)$.
- $D$ is bounded and strictly positive and $\Sigma_a$ is bounded on the domain $\Omega$ (Eq.(20)).
- $D$ is a piecewise differentiable function.

Under these assumptions, it follows that $\tilde{Q} \in L^2(\Omega)$. To derive the K-H equation, we multiply Eq.(8) by $G$ and Eq.(10) by $u$. Eqs.(53)(54) are obtained.

$$-G\nabla \cdot (D\nabla u) + G\Sigma_a u = G\tilde{Q} \quad (53)$$
$$-u\nabla \cdot (D\nabla G) + u\Sigma_a G = u\delta(\vec{x} - \vec{y}) \quad (54)$$

By subtracting Eq.(54) from Eq.(53), rearranging the terms and integrating over the domain $\Omega$, we obtain Eq.(55).

$$\int_\Omega [-G\nabla \cdot (D\nabla u) + u\nabla \cdot (D\nabla G)] \, d\Omega = \int_\Omega G\tilde{Q} \, d\Omega - \int_\Omega u\delta(\vec{x} - \vec{y}) \, d\Omega \quad (55)$$

By applying Green's Second Identity, the interior integrals can be transformed into boundary integrals over $\Gamma$, where $\vec{n}$ denotes the outward-pointing normal vector (Eqs.(56)(57)). The application of this identity highlights that the boundary terms emerge directly from the integration by parts. In particular, due to the homogeneous Dirichlet BC, the boundary integral in Eq.(57) vanishes.

$$\int_\Omega G\nabla \cdot (D\nabla u) \, d\Omega = \int_\Gamma GD\nabla u \cdot \vec{n} \, d\Gamma - \int_\Omega \nabla G \cdot (D\nabla u) \, d\Omega \quad (56)$$

$$\int_\Omega u\nabla \cdot (D\nabla G) \, d\Omega = \int_\Gamma uD\nabla G \cdot \vec{n} \, d\Gamma - \int_\Omega \nabla u \cdot (D\nabla G) \, d\Omega \quad (57)$$

Substituting these results into Eq.(55) leads to the formulation of Eq.(58).

$$u(\vec{x}) = \int_\Gamma G(\vec{x}, \vec{y})D(\vec{y})\nabla u(\vec{y}) \cdot \vec{n}(\vec{y})d\vec{y} + \int_\Omega G(\vec{x}, \vec{y})\tilde{Q}(\vec{y})d\vec{y} \quad \vec{x} \in \Omega \quad (58)$$

Let us examine the role of the Green's function in Eq.(58). For a Green's function that enforces homogeneous Dirichlet BC, the condition is imposed on the *field point*, i.e., the function is constructed to vanish when the field point lies on the boundary (Eq.(47)). On the right-hand side of Eq.(58), the first term is a boundary integral. Specifically, the field point $\vec{x}$ is kept free in the interior of the domain, while the integration is performed over the source variable $\vec{y}$, which lies on the boundary. The vanishing condition applies to the case where the field point is on the boundary - not the source point. So even with reciprocity, $G(\vec{x}, \vec{y})$ is generally nonzero for $\vec{y} \in \partial\Omega$ when $\vec{x} \in \Omega$. To better understand this behavior, let us revisit the reciprocity property stated in Eq.(49). Since Dirichlet BC is imposed on the field point, for a fixed field point $\vec{x}_1$, the Green's function satisfies Eq.(59).

$$G(\vec{x}_2, \vec{x}_1) = 0 \quad \vec{x}_2 \in \partial\Omega \tag{59}$$

By reciprocity, fixing $\vec{x}_2$ and treating $\vec{x}_1$ as the variable instead, the Green's function satisfies Eq.(60).

$$G(\vec{x}_1, \vec{x}_2) = 0 \quad \vec{x}_1 \in \partial\Omega \tag{60}$$

In the boundary integral representation, the roles of $\vec{x}_1$ and $\vec{x}_2$ are not interchangeable in the sense of "evaluation". With respect to the derived form of the K-H equation, although the homogeneous Dirichlet BC ensures that the field $u$ vanishes on the boundary, the boundary integral does not vanish entirely. This is due to contributions from both the normal derivatives and the structure of the Green's function itself.

### 4.2. *Procedure for reconstructing neutron flux under arbitrary control rod configurations*

The proposed methodology offers a significant computational advantage by shifting the focus from scenario-specific interpolation to the data-driven reconstruction of a fundamental physical object, i.e., the Green's function. A key property of the Green's function is that it does not depend on the external forcing term. This universality allows it to serve as a common descriptor across multiple operational scenarios. Unlike interpolation-based methods, which require large synthetic datasets covering a wide range of conditions and rely on training with scenarios similar to those expected during testing, the presented approach requires only a compact training dataset. Theoretically, measurements from a single reference scenario are sufficient to reconstruct the Green's function; once learned, it can be reliably applied to predict the neutron flux distribution under any other operational condition, as long as the geometry and the diffusive properties of the domain remain unchanged. This enables a "small-data" alternative to traditional big-data strategies. Rather than extracting mappings from large datasets, the method focuses on reconstructing a single, physics-informed pattern that accurately represents the neutron behavior within the reactor. For each new scenario, the problem is reformulated with homogeneous Dirichlet BC using the lifting method, ensuring consistency and general applicability of the learned Green's function. The procedure for applying the modified K-H equation is divided into two main stages, i.e., the training and the testing.

**Training Phase**
(Learning the Green's function from collected measurements)
- Step 1: Define the lifting function, $w$
    - Construct a function that satisfies the prescribed non-homogeneous Dirichlet BC.
    - The objective is to extend the boundary data ($\phi$ on $\partial\Omega$) smoothly into the interior of the domain using a linear function that matches measurements.

- Step 2: Evaluate the Effective Source term
  - Starting from the known source term ($Q$), the additive contributions introduced by the lifting method (Eq.(51)) are incorporated, and the effective source term $\tilde{Q}(\vec{y})$ is defined.

- Step 3: Compute the gradient of the unknown at the boundary
  - Collect boundary measurements of the neutron flux gradient.
  - Estimate the gradient of the solution, $\nabla u = \nabla \phi - \nabla w$ on $\partial \Omega$.

- Step 4: Solve the K-H equation for the Green's function
  - Treat the Green's function as the unknown in Eq.(58), and reconstruct it using the quantities derived from measurements, i.e., the values of $u$ and $\nabla u$ on $\partial \Omega$ and a few values of $u$ within $\Omega$. It is essential to complement boundary measurements with a sufficient number of interior point measurements to adequately characterize $G$.

**Testing Phase**
(Using the learned Green's function to reconstruct the neutron flux)

- Solve the K-H equation for the neutron flux
  - Use the previously reconstructed Green's function.
  - Solve Eq.(58) to estimate the neutron flux at any interior location of interest.
  - No interior point measurements are required. To reconstruct the field, it is sufficient to use the applied source term along with boundary measurements of the neutron flux and its gradient.

## 5. PROOF OF EXISTENCE AND UNIQUENESS OF THE GREEN'S FUNCTION OBTAINED FROM INVERSE PROBLEM SOLUTION

### 5.1. Formulation the Tikhonov-regularized problem

The problem to be addressed can be formulated as follows: "Given the measurements of $\phi$ and $\nabla \phi$ on $\partial \Omega$, a limited number of interior measurements of $\phi$, the full knowledge of the diffusive medium properties ($D$ and $\Sigma_a$) and the control rod position ($\tilde{Q}$), is it possible to uniquely and stably reconstruct the Green's function $G(\vec{x}, \vec{y})$?". To solve this problem, we propose of inverting Eq.(58). The inversion of an integral representation where $G$ appears as part of the kernel is not a standard coercive BVP in $G$. The operator mapping $G \mapsto u$ is not a simple bilinear form $G$, and Lax-Milgram lemma cannot be directly applied to prove the existence and uniqueness of $G$. As a result, this inverse problem is generally ill-posed in the sense of Hadamard. Some of the key challenges that may arise include:

- **Non-uniqueness**: Multiple Green's functions may produce the same observable data.
- **Instability**: Small measurement errors can lead to large variations in the reconstructed $G$.
- **Lack of regularity**: The solution may not belong to a well-defined function space without additional regularization.
- **Limited observability**: Having only boundary and a few interior measurements of $u$ is rarely sufficient to pin down the full function $G(\vec{x}, \vec{y})$, which is two-point-dependent.

To address issues of uniqueness and stability, additional integral constraints or model assumptions on $G$, complemented by sufficiently dense interior data and a regularization scheme, are required. In this work, we enforce the reciprocity property of the reconstructed Green's function (Section 3.5) and adopt a Tikhonov-type regularization [8], which reformulates the inverse problem as a minimization problem. The objective of this approach is to significantly reduce the set of admissible solutions and to yield a stable and well-defined optimal solution. In this Section, we provide a standard argument (often referred to as the *direct method in the calculus of variations* [9]) demonstrating that Tikhonov-regularized formulation of the problem satisfies existence, (potential) uniqueness, and stability. It is important to stress that this regularization does not make the problem well-posed in the classical Hadamard sense. By adding a Tikhonov regularization term and ensuring a sufficient number of measurements to adequately constrain $G$, a stable approximate solution can be obtained. This solution will exhibit continuous dependence on the input data, despite the presence of inherent measurement noise.

The *forward operator* $\mathcal{F}: X \to Y$ is defined in Eq.(61). It is a Fredholm integral operator of the first kind that maps the unknown Green's function - treated as an element of a suitable function space $X$ - into the data space $Y$, a normed space representing the measured data.

$$\mathcal{F}: G \mapsto u(\vec{x}, G) = \int_{\Gamma} G(\vec{x}, \vec{y}) D(\vec{y}) \nabla u(\vec{y}) \cdot \vec{n}(\vec{y}) d\vec{y} + \int_{\Omega} G(\vec{x}, \vec{y}) \tilde{Q}(\vec{y}) d\vec{y} \tag{61}$$

where $u(\vec{x}, G)$ is computed using the K-H representation.

As for the selection of the function and data spaces:

- The function space $X$ is chosen to ensure the regularity of the Green's function, guaranteeing that both the function and its weak derivatives (with respect to both variables) are square-integrable. $G(\vec{x}, \vec{y})$ is required to satisfy homogeneous Dirichlet BC in the $\vec{x}$-variable (Eq.(62)).

$$G(\vec{x}, \vec{y}) = 0 \quad \vec{x} \in \partial\Omega \tag{62}$$

The appropriate functional setting is the subspace of $H^1(\Omega \times \Omega)$ consisting of functions that vanish on the boundary in the $\vec{x}$-variable. This choice is particularly relevant since many Tikhonov functionals employ an $H^1$-norm as the regularization term, leading to a coercive and weakly lower semicontinuous (w.l.s.c.) penalty. The definition of this space is reported in Eq.(63).

$$X \equiv H^1_{0,x}(\Omega \times \Omega) = \{G \in H^1(\Omega \times \Omega) : G(\cdot, \vec{y}) \in H^1_0(\vec{y}) \text{ for almost every } \vec{y} \in \Omega\} \tag{63}$$

A natural norm on this space is given by Eq.(64). It is equivalent to the standard $H^1$ norm on $\Omega \times \Omega$ when restricted to functions that vanish on the boundary in the $\vec{x}$-variable.

$$\|G\|_X^2 = \|G\|_{H^1(\Omega \times \Omega)}^2 = \int_{\Omega} \left( \|G(\cdot, \vec{y})\|_{L^2(\Omega)}^2 + \|\nabla_{\vec{x}} G(\cdot, \vec{y})\|_{L^2(\Omega)}^2 + \|\nabla_{\vec{y}} G(\cdot, \vec{y})\|_{L^2(\Omega)}^2 \right) d\vec{y} \tag{64}$$

- Sensor data are naturally represented as square-integrable functions over their measurement region. Accordingly, the data space $Y$ is defined as $Y = L^2(\Omega)$ interior readings and $Y = L^2(\partial\Omega)$ for boundary readings.

The functional to be minimized is defined in Eq.(65).

$$J(G) = \{\|\mathcal{F}(G) - u_{obs}\|_Y^2 + \alpha \mathcal{R}(G)\} \tag{65}$$

with

- **Data-misfit term:** $\|\mathcal{F}(G) - u_{obs}\|_Y^2$. It measures how well the forward model predicts the observed data. Squaring the norm produces a differentiable and strictly convex functional.
- **Regularization term:** $\alpha \mathcal{R}(G)$. $\mathcal{R}(G)$ is a regularization functional that enforces smoothness and incorporates physical constraints, $\alpha > 0$ is the regularization parameter balancing the fit to the data against the smoothness and symmetry constraints.

*5.2. Definition of a suitable Regularization term*

In this application, the regularization term comprises two components: one enforcing smoothness and the other ensuring the reciprocity of the Green's function for a homogeneous Dirichlet problem (Eq.(45)).

- *Smoothness Term*

This term penalizes both large values of $G$ (through $\|G\|_{L^2(\Omega \times \Omega)}^2$), and rapid spatial variations in $G$ by taking gradients with respect to both of its arguments $\vec{x}$ and $\vec{y}$ (through $\|\nabla G\|_{L^2(\Omega \times \Omega)}^2$), thereby discouraging large oscillations. This term represents the norm of $G$ in function space $X$ (Eq.(66)).

$$\|G\|_X^2 = \|G\|_{L^2(\Omega \times \Omega)}^2 + \|\nabla G\|_{L^2(\Omega \times \Omega)}^2 = \int_{\Omega \times \Omega} \left( |G(\vec{x}, \vec{y})|^2 + |\nabla_{\vec{x}} G(\vec{x}, \vec{y})|^2 + |\nabla_{\vec{y}} G(\vec{x}, \vec{y})|^2 \right) d\vec{x} d\vec{y} \qquad (66)$$

- *Reciprocity Penalty*

This term penalizes deviations from the expected reciprocity condition by quantifying the asymmetry in $G$. By enforcing the expected symmetry of the Green's function as a soft constraint, it reduces the solution space while still permitting small deviations to account for noise or model inaccuracies. The simplest choice is an $L^2$- norm of the difference $\left(G(\vec{x}, \vec{y}) - G(\vec{y}, \vec{x})\right)$ (Eq.(67)).

$$\|G - G^T\|_{L^2(\Omega \times \Omega)}^2 = \int_{\Omega \times \Omega} |G(\vec{x}, \vec{y}) - G(\vec{y}, \vec{x})|^2 \, d\vec{x} d\vec{y} \qquad (67)$$

The regularization term is defined in Eq.(68).

$$\mathcal{R}(G) = \|G\|_{L^2(\Omega \times \Omega)}^2 + \|\nabla G\|_{L^2(\Omega \times \Omega)}^2 + \gamma \|G - G^T\|_{L^2(\Omega \times \Omega)}^2 = \|G\|_X^2 + \gamma \|G - G^T\|_{L^2(\Omega \times \Omega)}^2 \qquad (68)$$

where $\gamma > 0$ is a parameter that balances the strength of the symmetry enforcement relative to the smoothness requirement. The complete Tikhonov-regularized formulation for recovering $G$ is reported in Eq.(69).

$$\min_{G \in X} J(G) = \|\mathcal{F}(G) - u_{obs}\|_Y^2 + \alpha \left( \|G\|_X^2 + \gamma \|G - G^T\|_{L^2(\Omega \times \Omega)}^2 \right) \qquad (69)$$

The recovered function $G$ corresponds to the Green's function for the specific BVP under consideration, i.e., the neutron diffusion operator, on the given domain with homogenous Dirichlet BC. It serves as a "fundamental solution" in the sense that it is derived from the K-H equation using the definition of the Green's function (Eq.(10)), while remaining consistent with the measured data.

### 5.3. Properties of the forward operator

- *Linearity of $\mathcal{F}(G)$*

    We need to demonstrate that the operator in Eq.(61) is linear with respect to $G$. Let us define the functions $f_1$ and $f_2$ (Eqs.(70)(71)). These functions are known and are obtained from available measurements, BCs, or model parameters. Substituting these expressions into Eq.(61) leads to the expression in Eq.(72).

$$f_1(\vec{y}) = D(\vec{y})\nabla u(\vec{y}) \cdot \vec{n}(\vec{y}) \quad \vec{y} \in \partial\Omega \tag{70}$$
$$f_2(\vec{y}) = \tilde{Q}(\vec{y}) \quad \vec{y} \in \Omega \tag{71}$$
$$u(\vec{x}) = \int_\Gamma G(\vec{x}, \vec{y})f_1(\vec{y})d\vec{y} + \int_\Omega G(\vec{x}, \vec{y})f_2(\vec{y})d\vec{y} \tag{72}$$

    To prove linearity, we need to show that Eq.(73) holds for any scalars $\alpha, \beta$ and any two kernels $G_1, G_2$.

$$\mathcal{F}(\alpha G_1 + \beta G_2) = \alpha \mathcal{F}(G_1) + \beta \mathcal{F}(G_2) \tag{73}$$

    By plugging $\alpha G_1 + \beta G_2$ into $\mathcal{F}$, we obtain Eq.(74).

$$\mathcal{F}(\alpha G_1 + \beta G_2)(\vec{x}) = \int_\Gamma \alpha G_1 f_1 \, d\Gamma + \int_\Gamma \beta G_2 f_1 \, d\Gamma + \int_\Omega \alpha G_1 f_2 \, d\Omega + \int_\Omega \beta G_2 f_2 \, d\Omega \tag{74}$$

    By rearranging the terms, we obtain Eq.(75) for every $\vec{x}$, which corresponds to the definition of linearity.

$$\begin{aligned}\mathcal{F}(\alpha G_1 + \beta G_2)(\vec{x}) &= \alpha\left(\int_\Gamma G_1 f_1 \, d\Gamma + \int_\Omega G_1 f_2 \, d\Omega\right) + \beta\left(\int_\Gamma G_2 f_1 \, d\Gamma + \int_\Omega G_2 f_2 \, d\Omega\right) = \\ &= \alpha\mathcal{F}(G_1)(\vec{x}) + \beta\mathcal{F}(G_2)(\vec{x})\end{aligned} \tag{75}$$

- *Continuity (Boundedness) of $\mathcal{F}(G)$*

    To demonstrate that the $\mathcal{F}$ is bounded, we must prove that there exists a constant $L_\mathcal{F} > 0$ such that the condition in Eq.(76) is satisfied.

$$\|\mathcal{F}(G)\|_Y = \|u\|_{L^2(\Omega)} \leq L_\mathcal{F} \|G\|_X \tag{76}$$

    With respect to the operator in Eq.(72), we apply the Cauchy-Schwarz inequality to the first integral for each fixed $\vec{x}$, leading to Eq.(77).

$$\left|\int_\Gamma G(\vec{x}, \vec{y})f_1(\vec{y}) \, d\Gamma\right| \leq \int_\Gamma |G(\vec{x}, \vec{y})||f_1(\vec{y})| \, d\Gamma \leq \|G(\vec{x}, \cdot)\|_{L^2(\Gamma)} \|f_1\|_{L^2(\Gamma)} \tag{77}$$

    Now let us calculate the $L^2$ norm over $\vec{x} \in \Omega$ (Eq.(78)).

$$\left\|\int_\Gamma G(\vec{x}, \vec{y})f_1(\vec{y}) \, d\Gamma\right\|^2_{L^2(\Omega)} = \int_\Omega \left|\int_\Gamma G(\vec{x}, \vec{y})f_1(\vec{y})d\vec{y}\right|^2 d\vec{x} \leq \int_\Omega \|G(\vec{x}, \cdot)\|^2_{L^2(\Gamma)} \|f_1\|^2_{L^2(\Gamma)} \, d\vec{x} \tag{78}$$

Since $\|f_1\|^2_{L^2(\Gamma)}$ is independent of $\vec{x}$, we can derive the expression in Eq.(79). A similar argument applies to the second integral over $\Omega$ in Eq.(72), as shown in Eq.(80).

$$\left\| \int_\Gamma G(\vec{x},\vec{y})f_1(\vec{y})d\vec{y} \right\|_{L^2(\Omega)} \leq \|f_1\|_{L^2(\Gamma)} \|G\|_{L^2(\Omega\times\Gamma)} \tag{79}$$

$$\left\| \int_\Omega G(\vec{x},\vec{y})f_2(\vec{y})d\vec{y} \right\|_{L^2(\Omega)} \leq \|f_2\|_{L^2(\Omega)} \|G\|_{L^2(\Omega\times\Omega)} \leq \|f_2\|_{L^2(\Omega)} \|G\|_X \tag{80}$$

By combining the estimates and applying the triangle inequality in $L^2(\Omega)$, we obtain Eq.(81).

$$\|u\|_{L^2(\Omega)} \leq \|f_1\|_{L^2(\Gamma)} \|G\|_{L^2(\Omega\times\Gamma)} + \|f_2\|_{L^2(\Omega)} \|G\|_X \tag{81}$$

By virtue of a trace embedding inequality, the $L^2(\Omega \times \Gamma)$ norm can be controlled by the norm in the space $X$ (Eq.(82)), where $C_{\text{trace}}$ denotes the operator-norm of the trace map $H^1(\Omega) \to L^2(\Gamma)$ (Eq.(83)).

$$\|G\|_{L^2(\Omega\times\Gamma)} \leq C_{\text{trace}} \|G\|_X \tag{82}$$

$$C_{\text{trace}} = \sup_{G \in X\setminus\{0\}} \frac{\|G\|_{L^2(\Omega\times\Gamma)}}{\|G\|_X} \tag{83}$$

We can conclude that the operator $\mathcal{F}(G)$ is continuous from the function space $X$ into the data space $Y$ (Eq.(84)). The quantities that determine the continuity constant $L_\mathcal{F}$ ($\|f_1\|_{L^2(\Gamma)}$, $C_{\text{trace}}$ and $\|f_2\|_{L^2(\Omega)}$, as shown in Eq.(85)) are given by known functions, i.e., $f_1$ and $f_2$, as well as by the geometric properties of the domain, including the shape of $\Omega$, the regularity of its boundary, and the dimension of the ambient space.

$$\|\mathcal{F}(G)\|_Y \leq \left( \|f_1\|_{L^2(\Gamma)} C_{\text{trace}} + \|f_2\|_{L^2(\Omega)} \right) \|G\|_X \equiv L_\mathcal{F} \|G\|_X \quad \forall G \in X \tag{84}$$

$$L_\mathcal{F} = \|f_1\|_{L^2(\Gamma)} C_{\text{trace}} + \|f_2\|_{L^2(\Omega)} \tag{85}$$

### 5.4. Demonstration of well-posedness of the Tikhonov-regularized minimization problem

The Tikhonov-regularized problem is said to be well-posed if (1) a minimizer exists, (2) the minimizer is unique, and (3) the minimizer depends continuously on the collected measurements.

- *Existence of a minimizer*

To demonstrate the existence of a solution for the inverse problem, the direct method in the calculus of variations is used.

  o <u>Boundedness Below</u>

Since both $\|\mathcal{F}(G) - u_{obs}\|_Y^2$ and $\alpha \mathcal{R}(G)$ are non-negative, the functional $J(G) \geq 0, \forall G \in X$. Therefore, the infimum exists and finite (Eq.(86)). Hence, we can define a sequence $\{G_n\} \subset X$ such that the limit in Eq.(87) is attained. This sequence is called a *minimizing sequence*.

$$m = \inf_{G \in X} J(G) \tag{86}$$

$$\lim_{n \to \infty} J(G_n) = m \tag{87}$$

- Coercivity

A functional $J(G)$ is coercive in the space $X$ if $J(G) \to \infty$ as $\|G\|_X \to \infty$. In our case, the regularization term contains $\|G\|_X^2$ (Eq.(66)). This term grows at least as fast as $\|G\|_X^2$ for large $\|G\|_X$. Even though the forward operator $\mathcal{F}$ is not coercive by itself, its contribution in the Tikhonov functional is nonnegative, i.e., it does not "cancel out" the growth provided by the regularization term. Because of the dominance of the regularization term for large $\|G\|_X$, Eq.(88) follows directly.

$$J(G) \geq \alpha \|G\|_X^2 \tag{88}$$

This inequality guarantees that, if $J(G_n)$ stays near the minimal value, then the norms $\|G_n\|_X$ cannot blow up. This implies that the minimizing sequence $\{G_n\}$ is bounded in $X$ (Eq.(89)).

$$\|G_n\|_X \leq M \quad \forall n \tag{89}$$

- Weak lower semicontinuity

A functional $J(G)$ is w.l.s.c. if for every sequence $\{G_n\} \subset X$ that converges weakly to some $G \in X$ ($G_n \rightharpoonup G$), the inequality in Eq.(90) holds.

$$J(G) \leq \liminf_{n \to \infty} J(G_n) \tag{90}$$

In Hilbert spaces, any convex and continuous functional is w.l.s.c. Let us consider the regularization term in Eq.(68). Both components are squared norms. The mapping $G \mapsto \|G\|_X^2$ is convex and continuous on $X$. The term $\|G - G^T\|_{L^2(\Omega \times \Omega)}^2$ is a norm squared of a linear transformation of $G$, and it is convex and continuous as well. Then, if $G_n \rightharpoonup G$ in $X$, Eqs.(91)(92) hold.

$$\|G\|_X^2 \leq \liminf_{n \to \infty} \|G_n\|_X^2 \tag{91}$$

$$\|G - G^T\|_{L^2}^2 \leq \liminf_{n \to \infty} \|G_n - G_n^T\|_{L^2}^2 \tag{92}$$

Multiplying by the positive constants $\alpha$ and $\gamma$, we obtain Eq.(93).

$$\alpha \mathcal{R}(G) \leq \liminf_{n \to \infty} \alpha \left[ \|G_n\|_X^2 + \gamma \|G_n - G_n^T\|_{L^2}^2 \right] \tag{93}$$

As for the data-misfit term, the mapping $v \mapsto \|v - u_{obs}\|_Y^2$ is convex and continuous in the Hilbert space $L^2(\Omega)$, and then w.l.s.c. In Section 5.3, we demonstrated that the forward operator $\mathcal{F}: X \to Y$ is linear and continuous. As a consequence, the composite functional $G \mapsto \|\mathcal{F}(G) - u_{obs}\|_Y^2$ is a convex function of $G$ (the composition of a linear continuous operator with a convex function remains convex). As a result, this data-misfit term is w.l.s.c. Thus, if $G_n \rightharpoonup G$ in $X$, Eq.(94) holds. Since the sum of two w.l.s.c. functionals is w.l.s.c., it follows that Eq.(90) is satisfied.

$$\|\mathcal{F}(G) - u_{obs}\|_Y^2 \leq \liminf_{n \to \infty} \|\mathcal{F}(G_n) - u_{obs}\|_Y^2 \tag{94}$$

The Banach-Alaoglu theorem states that the closed unit ball in $X^*$ is compact in the weak$*$ topology. It guarantees the compactness needed to extract weakly convergent subsequences from bounded sequences in the dual space. Since $\{G_n\}$ is bounded, there exists a subsequence $\{G_{n_k}\}$ and an element $G^* \in X$ such that $G_{n_k} \rightharpoonup G^*$. Besides, given that the functional $J$ is w.l.s.c. and $\{G_{n_k}\}$ is a minimizing subsequence, Eq.(95) follows directly.

$$J(G^*) \leq \liminf_{k \to \infty} J(G_{n_k}) = m \tag{95}$$

By the definition of $m$ as the infimum, a solution to the regularized inverse problem exists, and $G^*$ is a minimizer of $J(G)$ (Eq.(96)).

$$J(G^*) = m \tag{96}$$

- *Uniqueness of the minimizer*

In Section 5.3, the linearity of the forward operator $\mathcal{F}$ was demonstrated. Since the norm squared of a linear function is a quadratic form, it follows that the data-misfit term is a convex function of $G$. The regularization term is strictly convex since both components are strictly convex. Therefore, $J(G)$ is strictly convex in $G$. Strict convexity implies that Eq.(97) is valid for any two distinct functions $G_1, G_2 \in X$ and any $0 < \lambda < 1$.

$$J(\lambda G_1 + (1-\lambda)G_2) < \lambda J(G_1) + (1-\lambda)J(G_2) \tag{97}$$

Let us suppose there were two distinct minimizers, i.e., $G_1^*$ and $G_2^*$. Then, using the property of strict convexity, we obtain Eq.(98).

$$J(\lambda G_1^* + (1-\lambda)G_2^*) < \lambda J(G_1^*) + (1-\lambda)J(G_2^*) \tag{98}$$

Since both $G_1^*$ and $G_2^*$ are minimizers, $J(G_1^*) = J(G_2^*) = m$. By substituting it in Eq.(98), we obtain Eq.(99).

$$J(\lambda G_1^* + (1-\lambda)G_2^*) < m \tag{99}$$

This inequality contradicts the minimality of both $G_1^*$ and $G_2^*$. Hence, the minimizer must be unique.

- *Stability with respect to the data*

In the previous Sections, we have shown that for any fixed data $u_{obs} \in Y$, the functional $J(G; u_{obs})$ is strictly convex, coercive and w.l.s.c. Hence, it has a unique minimizer $G(u_{obs})$ in the function space $X$. The ultimate objective of this work is to develop a methodology for real-time monitoring. Accordingly, the Green's function will be estimated using measurement data, which are inherently subject to experimental uncertainties. Let us consider a sequence of perturbed data $\{u_n\} \subset Y$ with $u_n \to u_{obs}$ in $Y$. For each $n$, let $G_n = G(u_n)$ denote the unique minimizer corresponding to the perturbed data, and $G^* = G(u_{obs})$ the unique minimizer corresponding to the unperturbed data. To demonstrate the stability of the solution with respect to the data, we need to show that the sequence of minimizers $\{G_n\}$ converges to $G^*$. For each $n$, we define $J(G; u_n)$ (Eq.(100)).

$$J(G; u_n) = \|\mathcal{F}(G) - u_n\|_Y^2 + \alpha \mathcal{R}(G) \tag{100}$$

By definition of $G_n$ and $G^*$, we have the inequalities in Eqs.(101)(102).

$$J(G_n; u_n) \leq J(G^*; u_n) \quad \forall n \tag{101}$$
$$J(G^*; u_{obs}) \leq J(G_n; u_{obs}) \tag{102}$$

Since the forward operator is linear and bounded, the mapping $G \mapsto \|\mathcal{F}(G) - u\|_Y^2$ is continuous with respect to $u$. In particular, for any fixed $G$, small perturbation in $u$ yield small changes in the data-misfit term (Eq.(103)).

$$\left| \|\mathcal{F}(G) - u_n\|_Y^2 - \|\mathcal{F}(G) - u_{obs}\|_Y^2 \right| \to 0 \text{ as } u_n \to u_{obs} \tag{103}$$

Since the regularization term does not depend on the data, this continuity means that the difference in the functional $|J(G_n; u_{obs}) - J(G_n; u_n)|$ can be bounded by a constant times $\|u_n - u_{obs}\|_Y^2$. Combining this observation with Eq.(102) yields the inequality in Eq.(104) up to an error that vanishes with $\|u_n - u_{obs}\|_Y^2$.

$$J(G_n; u_n) \approx J(G_n; u_{obs}) \geq J(G^*; u_{obs}) \tag{104}$$

Because the functional $J$ is strictly convex and coercive, the difference in the functional values between any two elements of $X$ provides control on the difference of the elements. Under strong convexity, it can be shown that there exists a constant $C > 0$ such that Eq.(105) holds.

$$\|G_n - G^*\|_X \leq C |J(G_n; u_{obs}) - J(G^*; u_{obs})|^{1/2} \tag{105}$$

Since the difference in the functional values due to data perturbation is small (by the continuity of the data-misfit term) and the minimizers are unique, Eq.(106) follows directly.

$$\|G_n - G^*\|_X \to 0 \text{ as } \|u_n - u_{obs}\|_Y^2 \to 0 \tag{106}$$

This shows that the minimizer depends continuously on the data, and the solution $G$ is stable with respect to small perturbations in $u_{obs}$.

## 6. EXTENSION OF WELL-POSEDNESS TO NONLINEAR SOURCE TERM SCENARIOS

In the analysis presented thus far, the source term has been assumed to be constant. Accounting for the dependence of local fission and absorption reaction rates on the neutron flux introduces a nonlinear coupling in the actual source term. Specifically, a more accurate expression for the source term, to be used in the algorithm, can be derived by combining Eqs.(6)(13)(51), and is presented in Eq.(107).

$$\tilde{Q} = \tilde{Q}(u) = (\nu \Sigma_f - \Sigma_a^{CR})(u + w) + \nabla \cdot (D \nabla w) - \Sigma_a w \tag{107}$$

A well-known technique for solving these problems consists of employing a lagged iterative approach [2]. The problem is recast as a Picard iteration, i.e., at each iteration, the source given by the previous flux estimate is treated as fixed, and the resulting linear diffusion equation is solved (Eq.(108)).

$$-\nabla(D \nabla u^n) + \Sigma_a u^n = \tilde{Q}(u^{n-1}) \tag{108}$$

where $u^n$ is the neutron flux approximation after iteration $n$, $\tilde{Q}(u^{n-1})$ is the source term computed using Eq.(107) and the flux estimate at the previous iteration ($u^{n-1}$).

In the context of solving the data-driven minimization problem, the nonlinearity of the source term does not pose an issue, provided that $\tilde{Q}$ acts as an admissible source, i.e., its values are consistent with and supported by the available measurements. However, employing a source term in the K-H equation that differs from the true physical source leads to an inaccurate forward operator $\mathcal{F}$. As a result, the predicted solution deviates from the true system behavior, and although the reconstructed Green's function may still minimize the cost functional, it will not accurately capture the underlying physical system. To establish the well-posedness of the neutron flux reconstruction problem with a nonlinear source term, the problem was reformulated as a sequence of linear problems with imposed sources (Eq.(109)). Based on the results presented in Section 5, if the convergence of the resulting sequence of linear problems can be demonstrated, then the well-posedness of the nonlinear formulation follows by induction.

$$
\begin{aligned}
u_{obs} &= \int_{\Gamma} G^1 D (\nabla u_{obs} \cdot \vec{n}) d\Gamma + \int_{\Omega} G^1 \tilde{Q}(u^0) d\Omega \\
u_{obs} &= \int_{\Gamma} G^2 D (\nabla u_{obs} \cdot \vec{n}) d\Gamma + \int_{\Omega} G^2 \tilde{Q}(u^1) d\Omega \\
u_{obs} &= \int_{\Gamma} G^n D (\nabla u_{obs} \cdot \vec{n}) d\Gamma + \int_{\Omega} G^n \tilde{Q}(u^{n-1}) d\Omega
\end{aligned}
\quad (109)
$$

where $u_{obs}$ represents the in-core measurements of the neutron flux and $\nabla u_{obs}$ represents the ex-core measurements of the gradient of the neutron flux (these values are fixed throughout the iterations), $u^0$ denotes the initial guess for the neutron flux distribution, and $G^n$ is the Green's function evaluated as the solution of the minimization problem at iteration $n$.

Table 2. Iterative procedure to estimate the Green's function when a nonlinear source term is applied.

| Iteration | Input data | Inverse problem to solve | Output data |
|---|---|---|---|
| 1 | $u^0$ | $u_{obs} = \mathcal{F}(G^1, \nabla u_{obs}, \tilde{Q}(u^0))$ | $G^1$ |
| 2 | $u^1 = u(G^1)$ | $u_{obs} = \mathcal{F}(G^2, \nabla u_{obs}, \tilde{Q}(u^1))$ | $G^2$ |
| $n$ | $u^{n-1} = u(G^{n-1})$ | $u_{obs} = \mathcal{F}(G^n, \nabla u_{obs}, \tilde{Q}(u^{n-1}))$ | $G^n$ |

The steps to solve the nonlinear problem are outlined in Table 2. The procedure begins with an initial guess for the neutron flux distribution ($u^0$). Based on this estimate, the corresponding source term ($\tilde{Q}(u^0)$) is computed and the minimization problem is solved. The resulting Green's function ($G^1$) is then used to reconstruct the neutron flux distribution throughout the domain, yielding an updated estimate $u^1$. This new flux is used to recalculate the source term, now scaled according to the improved flux ($\tilde{Q}(u^1)$). The process is repeated until convergence is achieved. Throughout the iterations, in addition to the properties of the diffusive medium, the values of the measured quantities ($u_{obs}$ and $\nabla u_{obs}$) and the Dirichlet BC used to compute the lifting function ($w$) remain unchanged. The only quantities that vary—and for which convergence is to be demonstrated—are the Green's function and the source term. The convergence of the iterative reconstruction scheme can be established by verifying the hypotheses of the Banach fixed-point theorem. First, the iterative process defined in Eq.(109) is reformulated as an operator $T: G^n \mapsto G^{n+1}$, which is composed of three successive sub-mappings (Eq.(110)).

$$G^{n+1} = T(G^n) = K^{-1}(S(\mathcal{F}(G^n))) \tag{110}$$

At iteration $n$, the neutron flux is calculated using the operator $\mathcal{F}: G^n \mapsto u^n$ defined in Eq.(61). Specifically, $u^n$ is reconstructed given $G^n$, the source $\tilde{Q}^{n-1}$ at previous iteration and the measured boundary gradient of the neutron flux imposed as input (Eq.(111)). The latter is given by $f_{1,obs}(\vec{y}) = D(\vec{y})\nabla u_{obs}(\vec{y}) \cdot \vec{n}(\vec{y})$ (Eq.(70)), which is treated as fixed.

$$u^n(\vec{x}) = \mathcal{F}(G^n) = \int_\Gamma G^n(\vec{x},\vec{y})f_{1,obs}(\vec{y})d\vec{y} + \int_\Omega G^n(\vec{x},\vec{y})\tilde{Q}^{n-1}(\vec{y})d\vec{y} \quad \forall \vec{x} \in \Omega \tag{111}$$

Given the calculated neutron flux, the source term is updated through the mapping $\mathcal{S}: u^n \mapsto \tilde{Q}^n$ (Eq.(112)).

$$\tilde{Q}^n = \tilde{Q}(u^n) = (\nu\Sigma_f - \Sigma_a^{CR})(u^n + w) + \nabla \cdot (D\nabla w) - \Sigma_a w \tag{112}$$

Finally, given the source term $\tilde{Q}^n$ as input, and with the in-core neutron flux measurements and the ex-core gradient measurements held fixed ($u_{obs}(\vec{x})$ and $f_{1,obs}(\vec{y})$), the operator $K^{-1}$ maps the source back to the Green's function ($K^{-1}: \tilde{Q}^n \mapsto G^{n+1}$), as defined in Eq.(113).

$$u_{obs}(\vec{x}) = \int_\Gamma G^{n+1}(\vec{x},\vec{y})f_{1,obs}(\vec{y})d\vec{y} + \int_\Omega G^{n+1}(\vec{x},\vec{y})\tilde{Q}^n(\vec{y})d\vec{y} \quad \forall \vec{x} \in \Omega \tag{113}$$

To establish the convergence of the iterative process, the Banach fixed-point theorem is applied. Since $X$ is a complete metric space, it is sufficient to show that the operator $T: X \mapsto X$ is a contraction with a constant $0 \leq L_T < 1$ such that $\|T(G_1) - T(G_2)\| \leq L_T\|G_1 - G_2\|$ for all $G_1, G_2 \in X$. Under this condition, the following results hold:

- A unique fixed point $G^* \in X$ exists (existence and uniqueness).
- The sequence of iterates converges, i.e., for any initial guess $G^0 \in X$, the iteration $G^{n+1} = T(G^n)$ converges to $G^*$.
- An error bound can be established: the convergence is at least geometric, with a rate $L_T$ (Eq.(114)).

$$\|G^{n+1} - G^*\|_X \leq \left(\frac{L_T^n}{1 - L_T}\right)\|G^1 - G^0\|_X \tag{114}$$

To verify that the operator $T$ is a contraction, it is necessary to prove that the mappings $\mathcal{F}$, $\mathcal{S}$ and $K^{-1}$ are each Lipschitz continuous.

- *Lipschitz continuity of $\mathcal{F}$ map*

In Section 5.3, the linearity and the boundedness of the forward operator $\mathcal{F}(G)$ were demonstrated (Eqs. (75)(84)). As $\mathcal{F}$ is a linear map between Banach spaces, its boundedness directly implies a global Lipschitz estimate (Eqs.(115)(116)).

$$\|\mathcal{F}(G_1) - \mathcal{F}(G_2)\|_Y = \|\mathcal{F}(G_1 - G_2)\|_Y \leq L_\mathcal{F}\|G_1 - G_2\|_X \quad \forall G_1, G_2 \in X \tag{115}$$

$$L_\mathcal{F} = \|f_1\|_{L^2(\Gamma)}C_{\text{trace}} + \|f_2\|_{L^2(\Omega)} \tag{116}$$

- *Lipschitz continuity of $\mathcal{S}$ map*

Since $\nabla \cdot (D\nabla w)$ and $\Sigma_a w$ are independent of $u$, and $\left(\nu\Sigma_f - \Sigma_a^{CR}\right) \in L^\infty(\Omega)$, we get Eq.(117). Thus, the $\mathcal{S}$ map is Lipschitz continuous from $Y$ into $L^2(\Omega)$ with constant $L_\mathcal{S}$ (Eq.(118)).

$$\left\|\tilde{Q}(u_1) - \tilde{Q}(u_2)\right\|_{L^2(\Omega)} \leq \left\|\nu\Sigma_f - \Sigma_a^{CR}\right\|_{L^\infty} \|u_1 - u_1\|_Y \equiv L_\mathcal{S} \|u_1 - u_1\|_Y \quad \forall u_1, u_2 \in Y \tag{117}$$

$$L_\mathcal{S} = \left\|\nu\Sigma_f - \Sigma_a^{CR}\right\|_{L^\infty} \tag{118}$$

- *Lipschitz continuity of $K^{-1}$ map*

The well-posedness of the Green's function reconstruction problem, formulated with boundary data and a prescribed source term, was established in Section 5.4. The functional minimized within the iterative scheme described in this section is defined in Eq.(119).

$$\min_{G^{n+1}\in X} J(G^{n+1}) = \left\|\left(\int_\Gamma G^{n+1}(\vec{x}, \vec{y}) D(\vec{y}) \nabla u_{obs}(\vec{y}) \cdot \vec{n}(\vec{y}) d\vec{y} + \int_\Omega G^{n+1}(\vec{x}, \vec{y}) \tilde{Q}^n(\vec{y}) d\vec{y}\right) - u_{obs}\right\|_Y^2 + \\ + \alpha \left(\|G^{n+1}\|_X^2 + \gamma \|G^{n+1} - (G^{n+1})^T\|_{L^2(\Omega\times\Omega)}^2\right) \tag{119}$$

We need to establish how sensitive the solution $G^{n+1}$ is to changes in the source term $\tilde{Q}^n$, while the observed data remains fixed. As shown in Eq.(119), $\tilde{Q}^n$ appears implicitly through the forward model. To derive the Lipschitz constant for the Green's function corresponding to the minimization problem, we need to consider the properties of the operator $K$ and the regularization used. The condition number of an operator $K$ is a measure of how sensitive the solution of a system involving $K$ is to changes in the input, i.e., the source term. For self-adjoint operators, the singular values of $K$ correspond to the absolute values of the eigenvalues of $K$. Hence, we can relate the condition number $\mathcal{K}(K)$ to the eigenvalues (Eq.(120)).

$$\mathcal{K}(K) = \frac{|\lambda_{\max}(K)|}{|\lambda_{\min}(K)|} \tag{120}$$

where $\lambda_{\max}(K)$ and $\lambda_{\min}(K)$ are the largest and the smallest eigenvalues of operator $K$, respectively. To this aim, the smallest eigenvalue $\lambda_{\min}(K)$ plays a crucial role in the stability of operator $K^{-1}$, i.e., a small $\lambda_{\min}(K)$ indicates that $K$ is ill-conditioned, meaning small changes in the data can lead to large changes in the solution.

The regularization parameters $\alpha$ and $\gamma$ in Tikhonov regularization helps shift the eigenvalues of $K$, and thus improves the conditioning of the inverse problem. Specifically, the smallest eigenvalue of $K$ gets shifted by the regularization parameters, which results in the regularized inverse operator being bounded. The Lipschitz constant for operator $K^{-1}$ ($L_{K^{-1}}$) is reported in Eq.(121).

$$L_{K^{-1}} = \frac{1}{\lambda_{\min}(K) + \alpha + \gamma} \tag{121}$$

By combining the Lipschitz continuity estimates for the three constituent mappings and applying Eq.(110), we obtain the inequality expressed in Eq.(122). The Lipschitz constant for the mapping $T$ is defined in Eq.(123).

$$\|G^{n+1} - G^n\|_X \leq L_{K^{-1}} \|\tilde{Q}^n - \tilde{Q}^{n-1}\|_{L^2(\Omega)} \leq L_{K^{-1}} L_S \|u^n - u^{n-1}\|_Y \leq L_{K^{-1}} L_S L_{\mathcal{F}} \|G^n - G^{n-1}\|_X \tag{122}$$

$$L_T = L_{K^{-1}} L_S L_{\mathcal{F}} = \frac{\left(\|f_1\|_{L^2(\Gamma)} C_{\text{trace}} + \|f_2\|_{L^2(\Omega)}\right) \|\nu \Sigma_f - \Sigma_a^{CR}\|_{L^\infty}}{\lambda_{\min}(K) + \alpha + \gamma} \tag{123}$$

If $L_T < 1$, then $T$ is a contraction on space $X$. By the Banach fixed-point theorem, $T$ admits a unique fixed-point ($G^* = T(G^*)$), and the sequence of iterates $G^n$ converge geometrically to $G^*$ (Eq.(114)). In this way, it can be demonstrated that the nonlinear problem can be reformulated as a convergent sequence of linear problems with imposed sources, for which well-posedness has already been established. Consequently, even when the source term is proportional to the local neutron flux, the solution to the optimization problem exists, is unique, and depends continuously on the measured data.

## 7. CONCLUSION

This work establishes the theoretical foundation for a novel methodology to reconstruct the neutron flux distribution in a nuclear reactor using a limited set of real-time, non-invasive measurements. The main idea is to approximate the Green's function—treated as a fundamental, physics-informed pattern governing neutron diffusion process—by training on a compact dataset associated with a single control rod configuration. Unlike interpolation-based approaches that require extensive training datasets, this method offers a "small-data" alternative. Rather than extracting mappings from datasets spanning diverse operational conditions, the proposed approach reconstructs a single Green's function tailored to a specific BVP (i.e., the neutron diffusion operator, on the given domain with homogenous Dirichlet BC) and consistent with the measured data. Through the lifting method, the BCs for any operational scenario can be mapped to a reference homogeneous Dirichlet scenario, ensuring the applicability of the learned Green's function across different control rod configurations.

As a major outcome, the well-posedness of the inverse problem is demonstrated, guaranteeing the existence and uniqueness of the Green's function and validating its use for neutron flux reconstruction. Regularity and symmetry properties are derived and integrated into the optimization algorithm, ensuring robustness against measurement noise and modeling uncertainties. The framework is organized in two main phases: a training phase, where the Green's function is inferred from reference data, and a testing phase, where it is used to reconstruct flux distributions under new operating scenarios. The methodology was first validated in the linear case with a constant source term and subsequently extended to the nonlinear case, where the source term is proportional to the local neutron flux—demonstrating its applicability under realistic reactor conditions.

### FUNDING

This work was supported by the U.S. Department of Energy, Office of Nuclear Energy. The submitted manuscript was created by UChicago Argonne, LLC, operator of Argonne National Laboratory. Argonne, a DOE Office of Science laboratory, is operated under contract DE-AC02-06CH11357.

access to these results of federally sponsored research in accordance with the DOE Public Access Plan at http://energy.gov/downloads/doe-public-accessplan (accessed on May 5, 2025).

## ACKNOWLEDGEMENTS

The author would like to express sincere gratitude to H. Wang (Argonne National Laboratory) for his valuable feedback and constructive suggestions, which significantly enhanced the quality and clarity of this paper.

## CONFLICTS OF INTEREST

The author declare that the research was conducted in the absence of any commercial or financial relationships that could be construed as a potential conflict of interest.

## REFERENCE

[1] A. Birri, K.C. Goetz, D.C. Sweeney, N.D. Bull Ezell, "Towards Realistic and High Fidelity Models for Nuclear Reactor Power Synthesis Simulation with Self-Powered Neutron Detectors", ORNL/TM-2023/3009, Oak Ridge National Laboratory, Oak Ridge, TN (2023).
[2] J.J. Duderstadt, L.J. Hamilton, "Nuclear Reactor Analysis", John Wiley and Sons, 1st Ed., New York, NY (1976).
[3] A.M. Quarteroni, A. Valli, "Numerical Approximation of Partial Differential Equations", Springer, 1st Ed., 2nd printing (2008).
[4] B. Kochunas, X. Huan, "Digital Twin Concepts with Uncertainty for Nuclear Power Applications", Energies, 14(14) (2021).
[5] J.A. Kulesza et al., "MCNP® Code Version 6.3.0 Theory & User Manual", LA-UR-22-30006, Rev. 1, Los Alamos National Laboratory, Los Alamos, NM (2022).
[6] H. Wang, R. Ponciroli, R.B. Vilim, V. Theos, S. Chatzidakis, "A Data-driven approach to Core Power distribution reconstruction in a Nuclear Reactor", ANL/NSE-24/46, Argonne National Laboratory, Lemont, Illinois (2024).
[7] J.D. Jackson, "Classical Electrodynamics", John Wiley and Sons, 3rd Ed., New York, NY (1999).
[8] G.H. Golub, P.C. Hansen, D.P. O'Leary, "Tikhonov Regularization and Total Least Squares", SIAM Journal on Matrix Analysis and Applications, 21(1), 184-194 (1999).
[9] B. Dacorogna, "Direct Methods in the Calculus of Variations", Springer-Verlag, 2nd Ed., New York, NY, (2007).